\documentclass[12pt]{amsart}

\usepackage{amsfonts,amsmath,amssymb}

\title{Homotopy theory of Hopf Galois extensions}
\author{Christian Kassel}
\address{Christian Kassel: Institut de Recherche Math\'{e}matique Avanc\'{e}e,
CNRS-Universit\'{e} Louis Pasteur,
7 rue Ren\'{e} Descartes, 67084 Strasbourg, France}
\email{kassel@math.u-strasbg.fr}
\author{Hans-J\"urgen Schneider}
\address{Hans-J\"urgen Schneider:
Mathematisches Institut, Universit\"at M\"un\-chen, Theresienstr. 39,
D-80333 Munich, Germany}
\email{Hans-Juergen.Schneider@mathematik.uni-muenchen.de}

\newtheorem{Lem}{Lemma}[section]
\newtheorem{Prop}[Lem]{Proposition}
\newtheorem{Cor}[Lem]{Corollary}
\newtheorem{Thm}[Lem]{Theorem}

\theoremstyle{definition}
\newtheorem{Def}[Lem]{Definition}
\newtheorem{Rem}[Lem]{Remark}
\newtheorem{Rems}[Lem]{Remarks}
\newtheorem{Expl}[Lem]{Example}
\newtheorem{Expls}[Lem]{Examples}

\newcommand\pf{\begin{proof}}
\newcommand\epf{\end{proof}}

\newcommand\Gal{\operatorname{Gal}}
\newcommand\CGal{\operatorname{CGal}}
\newcommand\Pic{\operatorname{Pic}}

\newcommand\Hom{\operatorname{Hom}}

\newcommand\Coker{\operatorname{Coker}}

\newcommand\id{{\operatorname{id}}}
\newcommand\can{{\operatorname{can}}}
\newcommand\co{{\operatorname{co}}}

\newcommand\cop{{\operatorname{cop}}}

\renewcommand\o{\otimes}

\newcommand\sw[1]{{}_{(#1)}}

\numberwithin{equation}{section}

\begin{document}

\maketitle

\vskip 30pt
\noindent
{\sc Abstract}.
{\it We introduce the concept of homotopy equivalence
for Hopf Galois extensions and make a systematic study of it.
As an application we determine all $H$-Galois extensions 
up to homotopy equivalence in the case when $H$ is a Drinfeld-Jimbo
quantum group.
}


\bigskip
\noindent
{\sc Titre} : 
{\bf Homotopie des extensions de Hopf Galois}

\bigskip
\noindent
{\sc R\'esum\'e}.
{\it Nous \'etudions le concept d'\'equivalence d'homotopie
pour les extensions $H$-galoisiennes o\`u $H$ d\'esigne une alg\`ebre de 
Hopf. Ceci nous permet de classifier les extensions $H$-galoisiennes
\`a homotopie pr\`es lorsque $H$ est un groupe quantique de Drinfeld-Jimbo.}

\bigskip
\noindent
{\sc Mots-cl\'es:}
extension galoisienne, alg\`ebre de Hopf, groupe quantique,
homotopie, g\'eom\'etrie non commutative, fibr\'e principal


\bigskip\bigskip

\section*{Introduction}

The purpose of this article is to make a systematic study of the concept
of homotopy equivalence introduced in the framework of Hopf Galois extensions
by the first-named author~\cite{K}.
As has been stressed many times (see e.g.~\cite{S}), Hopf Galois
extensions can be viewed as non-commutative analogues of principal
fibre bundles where the role of the structural group is played by a
Hopf algebra. It is therefore natural to adapt the concept of homotopy to them.

Hopf Galois extensions for a given Hopf algebra
and over a given algebra are difficult to classify up to isomorphism.
One of our motivations in this paper and in~\cite{K} was that it might be
easier to classify Hopf Galois extensions up to homotopy equivalence.
We show in this paper that it is indeed so in the case
when $H$ is a Drinfeld-Jimbo quantum group or some finite-dimensional variant.
More precisely, we prove that the homotopy classes of $H$-Galois extensions
for such a Hopf algebra $H$ are in bijection with the
homotopy classes of $k[G]$-Galois extensions,
where $G$ is the group of group-like elements of~$H$.

Certain $K$-theoretic elements naturally attached
to Hopf Galois extensions and recently investigated in connection with
non-commutative geometry (see the survey article~\cite{BH} and references
therein)
turn out to be homotopy invariant. This gives another reason to
consider homotopy equivalence in this algebraic framework.

In Section~1 we collect basic definitions and various results in the
literature in order to state functorial properties of the set
$\Gal_B(H)$ of isomorphism classes of faithfully flat $H$-Galois
extensions of an algebra~$B$. This includes change of scalars, change
of Hopf algebras, and twistings of Galois extensions.

Section~2 explains the concept of homotopy equivalence of Hopf
Galois extensions, which had been introduced in~\cite{K} for central
extensions (contrary to {\em loc.\ cit.}, we avoid here any reference to
\'etale morphisms).  Our most striking result on the set $\mathcal{H}_B(H)$
of homotopy classes of faithfully flat $H$-Galois
extensions of an algebra~$B$ is the following:
if $H = \oplus_{n\geq 0}\, H(n)$ is an $\mathbb{N}$-graded Hopf algebra,
then the inclusion $H(0) \subset H$ induces a bijection
$$\mathcal{H}_B(H) \cong \mathcal{H}_B(H(0)).$$

In Section~3 we make the connection between our definition of
homotopy equivalence and the construction of homotopy functors
in algebraic $K$-theory from which it is directly inspired.
When $B$ is a commutative ring satisfying certain conditions and
$H = k[G]$ is a group algebra, we relate $\Gal_B(H)$ and $\mathcal{H}_B(H)$
to the Picard group $\Pic(B)$ and some cohomology group of~$G$.

Finally, in Section~4 we show that, when $U_q(\mathfrak{g})$ is a
Drinfeld-Jimbo quantum group and $G$ is the group of its group-like elements,
then
$$\mathcal{H}_B(U_q(\mathfrak{g})) \cong \mathcal{H}_B(k[G]).$$
This follows from the results obtained in the previous sections and from
the fact that $U_q(\mathfrak{g})$ is a twist of a graded Hopf algebra
whose component of degree~0 is~$k[G]$
(for results related to the latter fact, see \cite{AS},~\cite{Di}).
\medskip

Throughout the paper, the ground ring is an arbitrary commutative ring~$k$.
Unadorned tensor product means tensor product over~$k$,
and algebras and coalgebras are defined over~$k$.
We denote by $U(R)$ the group of invertible elements of a ring~$R$,
by $G(H)$ the group of group-like elements of a Hopf algebra~$H$,
and by $C_n$ the cyclic group of order~$n$.

If $C$ is a coalgebra, and $V$ a right (resp.\ left) $C$-comodule,
we use the following version of the Sweedler notation:
$\Delta(c) = c\sw1 \o c\sw2$ for the comultiplication of $c \in C$,
and $v\sw0 \o v\sw1$ (resp.\ $v\sw{-1} \o v\sw0$) for the coaction of~$v
\in V$.

\section{Hopf Galois extensions}

Let $H$ be a Hopf algebra, and $A$ a right $H$-comodule algebra with
structure map
$\delta : A \to A \o H, \; a \mapsto \delta(a) = a\sw0 \o a\sw1$,
that is, $\delta$ is an algebra map and a right $H$-comodule structure.
The subalgebra $A^{\co H} \subset A$ of $H$-coinvariant elements is defined as
$$A^{\co H}= \{a \in A \mid a\sw0 \o  a\sw1 = a \o 1\}.$$
Let $B$ be an algebra. A (right) $H$-comodule algebra over~$B$
is a triple $(A,\delta, i)$,
where $A$ is a right $H$-comodule algebra with structure map~$\delta$
and $i : B \to A$ is an algebra map with $i(B) \subset A^{\co H}$.
For any $H$-comodule algebra $(A,\delta, i)$ over~$B$, the canonical or
Galois map
is defined~by
$$\can : A \o_B A \to A \o H,\; x\o y \mapsto xy\sw0 \o y\sw1.$$
An $H$-comodule algebra $(A,\delta, i)$, or simply $A$ over $B$, is an
{\em $H$-Galois extension of~$B$} if its Galois map is bijective
and if $i$ defines an isomorphism $i : B \to A^{\co H}$.
If $A$ is a right $H$-comodule algebra and $B$ is a subalgebra of $A$,
we will say that $B \subset A$ is an $H$-Galois extension if $B = A^{\co H}$
and if the Galois map of $A$ over~$B$ is bijective.

An $H$-comodule algebra $A$ is called $H$-{\em cleft}
if there exists a right $H$-colinear map $\gamma : H \to A$
that is invertible up to convolution.
If $A$ is $H$-cleft, then $A^{\co H} \subset A$ is $H$-Galois
(see~\cite[8.2.4]{M}).

A morphism
$\varphi : (A,\delta, i) \to (\widetilde{A},\widetilde{\delta},
\widetilde{i})$
of $H$-comodule algebras over $B$ is a right $H$-colinear algebra map
$\varphi : A \to \widetilde{A}$ with $\widetilde{i} = \varphi i$.
An isomorphism of $H$-comodule algebras is a bijective morphism.
Note that an $H$-colinear algebra isomorphism of $H$-Galois extensions of $B$
is an isomorphism of comodule algebras over~$B$.

\begin{Lem}\label{GaloisLem}
Let $H$ be a Hopf algebra, $B$ be an algebra, and $(A,\delta,i)$ and
$(\widetilde{A},\widetilde{\delta}, \widetilde{i})$ be $H$-comodule
algebras over~$B$
with bijective Galois maps.
\begin{enumerate}
\item \cite[4.2]{T}
Assume that $A$ is faithfully flat as a right or left $B$-module via~$i$.
Then $A$ is an $H$-Galois extension of $B$.
\item \cite[3.11(1)]{S}
Let $\varphi : A \to \widetilde{A}$ be a morphism of $H$-comodule algebras
over~$B$.
Assume that $\widetilde{A}$ is faithfully flat as a right $B$-module
via~$\widetilde{i}$.
Then $\varphi$ is an isomorphism, and both $A$ and $\widetilde{A}$
are $H$-Galois extensions of~$B$.
\end{enumerate}
\end{Lem}

To define homotopy (see Section 2) we have to extend the ground ring from $k$
to the polynomial algebra $k[t]$ in the indeterminate~$t$.

In general, let $\alpha : k \to R$ be a homomorphism of commutative rings.
Ground ring extension from $k$ to $R$ then means tensoring with $R$ over~$k$,
where $R$ is a module over $k$ via~$\alpha$.
If  $B$ is a $k$-algebra and $M$ a left $B$-module, then $R \o B$ is an
$R$-algebra, and $R \o M$ is a left $R \o B$-module in the natural way by
$$(r \o b)(s \o m) = rs \o bm$$
for all $r,s \in R$, $b \in B$, $m \in M$.
Recall that for any right $R\o B$-module~$X$,
$$X \o_{R \o B} (R \o M) \to X \o_B M, \; x \o r \o m \mapsto xr \o m,$$
is an isomorphism.
In particular, ground ring extension preserves flatness and faithful flatness.
Similarly, if $C$ is a $k$-coalgebra and $V$ a right $C$-comodule,
then $R \o C$ is an $R$-coalgebra and $R \o V$ is a right $R \o C$-comodule.
Note that a right $R \o C$-comodule structure on an $R$-module~$W$ is an
$R$-linear right $C$-comodule structure
$$W \to W \o_R (R \o C) \cong W \o C.$$
If $A$ is an $H$-comodule algebra over $B$, then we obtain by ground ring
extension
an $R \o H$-comodule algebra $R\o A$ over $R\o B$ with ground ring~$R$.

\begin{Prop}\label{baseextension}
Let $H$ be a Hopf algebra, $B$ an algebra, and $A$ an $H$-Galois extension
of~$B$
with ground ring~$k$. Assume that $A$ is left (resp.\ right) faithfully
flat over~$B$.
Let $\alpha : k \to R$ be a homomorphism of commutative rings.
Then $R \o A$ is an $R \o H$-Galois extension of~$R \o B$, and
$R \o A$ is left (resp.\ right) faithfully flat over~$R \o B$.
\end{Prop}

\pf
It is clear that the Galois map of $R \o A$ is bijective and that $R \o A$ is
left (resp.\ right) faithfully flat over~$R\o B$.
Hence the claim follows from Lemma \ref{GaloisLem}~(1).
\epf

The functorial behaviour of $H$-Galois extensions or $H$-comodule algebras
in the Hopf algebra $H$ as a variable is given by the cotensor product.
If $C$ is a coalgebra, and $V$ and $W$ are right and left $C$-comodules
with comodule structures $\delta_V : V \to V \o C$ and
$\delta_W : W \to C \o W$, then the cotensor product $V \Box_C W$ is the
kernel of
$$\delta_V \o \id_W - \id_V \o \delta_W : V \o W \to V \o C \o W.$$
Let $\varphi : K \to H$ be a Hopf algebra homomorphism, and $A$ a right
$H$-comodule
algebra over $B$ with algebra map~$i : B \to A$.
Assume that $K$ is flat as a $k$-module. Then $A \Box_H K \subset A \o K$
is a subalgebra and right $K$-subcomodule of~$A \o K$ with
componentwise multiplication and $K$-comodule structure~$\id_A \o \Delta$.
Here we view $K$ as a left $H$-comodule by
$$K \to H \o K, \;x \mapsto \varphi(x\sw1) \o x\sw2.$$
Note that the inclusion $A \Box_H K \subset A \o K$ defines an injective map
$$(A \Box_H K) \o K \to A \o K \o K$$
since $K$ is flat over~$k$.
Hence the $K$-comodule structure of $A \Box_H K$ is well defined.

Moreover, $A \Box_H K$ is a $K$-comodule algebra over $B$ by the algebra map
$B \to A \Box_H K, \; b \mapsto i(b) \o 1$.

\begin{Prop}\cite[3.11 (3)]{S}\label{cotensor}
Let $\varphi : K \to H$ be a Hopf algebra homomorphism
and $A$ an $H$-Galois extension of $B$.
Assume that $A$ is right faithfully flat over $B$ and that $K$ is flat
over~$k$.
Then $A \Box_H K$ is a $K$-Galois extension of $B$ and
right faithfully flat over~$B$.
\end{Prop}

Let $H$ be a Hopf algebra and $B$ an algebra over $k$. We denote
by~$\Gal_B(H/k)$,
or simply by~$\Gal_B(H)$, the set of isomorphism classes of $H$-Galois
extensions of $B$ that are right faithfully flat over~$B$.

Let $\alpha : R \to S$ be a homomorphism of commutative $k$-algebras.
For any right faithfully flat $R \o H$-Galois extension $A$ of $R \o B$
let $\alpha_*A = S \o_{R} A$ be the $S \o H$-Galois extension of $S \o B$
defined by ground ring extension, where we identify
$$S \o_{R} (R \o H) \cong S \o H \text{ and } S \o_{R} (R \o B) \cong S \o
B.$$

Let $\varphi : K \to H$ a Hopf algebra homomorphism. Assume that $K$ is
flat over~$k$.
For any right faithfully flat $H$-Galois extension $A$ of~$B$
let  $\varphi^*A = A \Box_H K$ be the
$K$-Galois extension of $B$ given by the cotensor product.

Then $\alpha_*$ and $\varphi^*$ define maps
\begin{align*}
\alpha_* : \Gal_{R \o B}(R \o H)& \to \Gal_{S \o B}(S \o H),\notag\\
\varphi^* : \Gal_B(H)& \to \Gal_B(K).\notag
\end{align*}

We collect the basic rules concerning $()^*$ and $()_*$.
In particular, $\Gal_{R \o B}(R \o H)$ is a covariant functor in
commutative $k$-algebras~$R$ and $\Gal_B(H)$ is a contravariant functor
in $k$-flat Hopf algebras~$H$.

\begin{Prop}\label{commute}
Let $H$ be a Hopf algebra over~$k$ and $A$ a right faithfully flat $H$-Galois
extension of~$B$.
Let $\alpha : k \to R$ and $\beta : R \to S$ be homomorphisms of
commutative $k$-algebras, $K$, $L$ be $k$-flat Hopf algebras,
and $\psi : L \to K$, $\varphi : K \to H$ be Hopf algebra homomorphisms.
Then
\enumerate
\item $(\alpha \beta)_*A \cong \alpha_* \beta_*A$.
\item $(\varphi \psi)^*A \cong \psi^* \varphi^*A$.
\item $\alpha_*\varphi^*A \cong (\id_R \o \varphi)^* \alpha_*A$.
\end{Prop}

\pf
(1) follows from the associativity of the tensor product.

To prove (2), we first define a right $L$-colinear algebra map
$$f : A \o L \to A \o K \o L, \; a \o x \mapsto a \o \psi(x\sw1) \o x\sw2.$$
The defining sequence
$$0 \to A \Box_HK \to A \o K \to A \o H \o K$$
of $\varphi^*A = A \Box_HK$ remains exact after tensoring with $L$ over~$k$
since $L$ is $k$-flat. Then we see that $f$ maps the subspace
$(\varphi \psi)^*A = A \Box_HL$ into~$(A \Box_HK) \o L$.
Finally, $f(A \Box_HL)$ is contained in 
$\psi^* \varphi^*A = (A \Box_H K)\Box_K L$
by the definition of the cotensor product of $A \Box_H K$ with~$L$.
Note that we have used flatness of $K$ over $k$ to define the $k$-comodule
structure of~$A \Box_H K$.

Thus we have defined a right $L$-colinear algebra map
$$A \Box_H L \to (A \Box_H K) \Box_K L$$
which is the identity on the subalgebra $B$ embedded in the first factor.
By Lemma \ref{GaloisLem}~(2) and Proposition \ref{cotensor},
this map is an isomorphism, and (2) is proved.

(3) By Propositions \ref{baseextension} and Proposition \ref{cotensor},
$$\alpha_* \varphi^*A = R \o (A \Box_H K) \text{ and }
(\id_R\o \varphi)^*\alpha_*A=(R \o A) \Box_{R \o H} (R \o K)$$
are both right faithfully flat $R \o K$-Galois extensions of~$R \o B$.
Since the map
$$R \o (A \Box_H K) \to (R \o A) \Box_{R \o H} (R \o K) $$
defined by
$r \o \sum_i a_i \o x_i \mapsto \sum_i r  \o a_i \o 1 \o x_i$
is well defined and a morphism of $R \o K$-Galois extension of~$R \o B$,
it is bijective by Lemma~\ref{GaloisLem}~(2).
\epf

Let $H$ be a Hopf algebra and $R$ a commutative algebra.
An $H$-Galois extension $(A,\delta,i)$ of  $R$ is called {\em central}
if $i(R)$ is contained in the center of~$A$.
We denote by~$\CGal_R(H/k)$, or simply by~$\CGal_R(H)$,
the set of isomorphism classes of faithfully flat central $H$-Galois
extensions of~$R$.

Faithfully flat $H$-Galois extensions of the ground ring~$k$ are also called
$H$-{\em Galois objects}.

\begin{Rem}
If $R \subset A$ is an $H$-Galois extension and if $R$ is central in~$A$,
then $A$ is an $R$-algebra and we can view $R \subset A$ as an $R \o H$-Galois
extension of~$R$ over the ground ring~$R$ with comodule algebra structure
$$A \to A \o H \cong A \o_R (R \o H).$$

Thus central Galois extensions of $R$ can be identified with Galois objects
and there is a functorial isomorphism
$$R \mapsto  \Gal_R(R \o H/R) \cong \CGal_R(H/k).$$
\end{Rem}

\medskip

Finally we briefly recall the notion of twisting of Hopf algebras and
comodule algebras \cite[Theorem 1.6]{D},~\cite[10.2.3]{KS}.
Let $H$ be a Hopf algebra and $\sigma : H \o H \to k$ be a (normalized)
{\em 2-cocycle}, that is,
\begin{equation}\label{cocycle}
\sigma(x\sw1,y\sw1)\, \sigma(x\sw2y\sw2,z)
= \sigma(y\sw1,z\sw1)\, \sigma(x,y\sw2z\sw2),
\end{equation}
for all $x,y,z \in H$, and $\sigma(x,1) = \varepsilon(x) = \sigma(1,x)$.

Assume that $\sigma$ is invertible with inverse $\sigma^{-1}$ with respect
to the convolution product.
Then the twisted Hopf algebra $H^{\sigma}$ is $H$ as a coalgebra with
the twisted multiplication
\begin{equation}\label{twisted}
x\cdot_{\sigma} y = \sigma(x\sw1,y\sw1)\, x\sw2y\sw2 \,
\sigma^{-1}(x\sw3,y\sw3).
\end{equation}

Let $A$ be a right $H$-comodule algebra. Then the twisted
$H^{\sigma}$-comodule algebra $A^{\sigma}$ is $A$ as an $H$-comodule
with the twisted multiplication
\begin{equation}\label{twisted2}
a \cdot_{\sigma} b = a\sw0 b\sw0 \, \sigma^{-1}(a\sw1,b\sw1).
\end{equation}

The convolution inverse $\sigma^{-1}$ of $\sigma$ is an invertible 2-cocycle
for~$H^{\sigma}$, and $(A^{\sigma})^{\sigma^{-1}} = A$ as a comodule
algebra over
$(H^{\sigma})^{\sigma^{-1}} = H$.

\begin{Prop}\label{twisting}
Let $H$ be a Hopf algebra, $\sigma : H \o H \to k$ an invertible 2-cocycle,
and $B$ an algebra.
\begin{enumerate}
\item If $A$ is an $H$-Galois extension of $B$, then $A^{\sigma}$ is an
$H^{\sigma}$-Galois extension of~$B$.
\item The map $\Gal_B(H) \to \Gal_B(H^{\sigma})$ given by $A \mapsto
A^{\sigma}$
is bijective.
\end{enumerate}
\end{Prop}

\pf
(1) is shown in~\cite[Section 4]{MS}, and (2) follows from (1)
since $(A^{\sigma})^{\sigma^{-1}} =A$.
\epf

 \section{Homotopy}

The homotopy properties of principal bundles in topology are based on the
following result. Let $G$ be a topological group, $B$ a topological space,
and $\xi$ a numerable principal $G$-bundle over~$B$.
If $f_i : B' \to B$, $i =0,1$, are homotopic maps, then the pull-back bundles
$f_0^*(\xi)$ and $f_1^*(\xi)$ are isomorphic \cite[Chapter 4, Theorem
9.9]{Hu}.

An $H$-Galois extension $B \subset A$ can be viewed as the algebraic analogue
of the function algebra of a principal $G$-bundle $p : X \to B$.
To study the algebraic analogue of homotopy for Hopf Galois extension
we first recall the definition of homotopic algebra maps used in algebraic
$K$-theory (see \cite[Section~3]{G},~\cite[Section~4]{Sw1}).

Let $k \subset k[t]$ be the polynomial algebra in the indeterminate~$t$.
For any $k$-module $V$, we denote the ground ring extension with respect
to $k \subset k[t]$ by~$V[t]$.
We define $k$-linear maps $V[i] : V[t] \to V$ for $i \in \{0,1\}$
by sending $vt^n$ to~$vi^n$. We will usually write simply $[i]$ instead
of~$V[i]$.
Note that these maps are algebra maps if $V$ is an algebra.

\begin{Def}
Let $\alpha, \beta : R \to S$ be homomorphisms of commutative $k$-algebras.
We say that $\alpha$ and $\beta$ are {\em homotopic} ($\alpha \sim \beta$)
if there exists a $k$-algebra homomorphism $\theta : R \to S[t]$ such that
$$[0] \theta = \alpha \text{ and } [1] \theta = \beta.$$
The map $\theta$ is called a {\em homotopy} between $\alpha$ and~$\beta$.
Let $\approx$ be the equivalence relation generated by~$\sim$.

The homomorphism $\alpha : R \to S$ is a {\em homotopy equivalence} if
there exists a ring homomorphism $\alpha' : S \to R$ with
$\alpha \alpha' \approx \id_{S}$ and $\alpha' \alpha \approx \id_R$.
\end{Def}

\begin{Rems}\label{homotopicmaps}
(1) The relation $\sim$ is reflexive and symmetric
(as one sees by the mapping $t \mapsto t-1$).
Thus $\alpha \approx \beta$ means that there is a sequence of algebra
homomorphisms
$\alpha_1, \dots, \alpha_n : R \to S$ with $\alpha_1 = \alpha$,
$\alpha_1 \sim \alpha_2, \dots, \alpha_{n-1} \sim \alpha_n$,
and~$\alpha_n = \beta$.

(2) Let $R = \oplus_{n\geq 0} R(n)$ be an $\mathbb{N}$-graded commutative
algebra.
Then the inclusion $R(0) \subset R$ is a homotopy equivalence of
commutative rings.
We recall the argument in~\cite{K}. Let $\iota : R(0) \to R$ be the inclusion,
and $\pi : R \to R(0)$ the projection given by the grading.
Then $\iota$ and $\pi$ are ring homomorphisms, and $\pi \iota = \id_{R(0)}$.
We define an additive map $\theta : R \to R[t]$ by
$\theta(x)= t^nx$ if $x \in R(n)$ ($n\geq 0$).
Since $R$ is $\mathbb{N}$-graded, $\theta$ is a ring homomorphism.
By construction, $[0]\theta = \iota \pi$ and $[1]\theta = \id_R$,
hence $\iota \pi \sim \id_R$.

(3) Here is a non-graded example of a homotopy equivalence.
Let $k$ be a field of characteristic $p >0$ and $k[G]$ the group algebra of a
finite abelian $p$-group $G$.
Write $G$ as a direct product of cyclic $p$-groups with generators
$g_i, 1 \leq i \leq r$. Then the algebra map
$$\theta : k[G] \to k[G][t] \text{ with }
\theta(g_i) = 1 + t(g_i-1), 1 \leq i \leq r,$$
is a homotopy with
$[0]\theta(g_i) =1$ and $[1]\theta(g_i) = g_i$ for all~$i$.
Hence the inclusion $k \subset k[G]$ is a homotopy equivalence.
\end{Rems}

\begin{Def}
Let $H$ be a Hopf algebra and $B$ an algebra.
Let $A_0$ and $A_1$ be right faithfully flat $H$-Galois extensions of~$B$.
We write $A_0 \sim A_1$ if there exists a right faithfully flat
$H[t]$-Galois extension $A$ of~$B[t]$ with ground ring~$k[t]$
such that $[i]_*A \cong A_i$ for  $i \in \{0,1\}$.

{\em Homotopy equivalence} of right faithfully flat $H$-Galois extensions
of~$B$ is the equivalence relation $\approx$ generated by~$\sim$.

Let $\mathcal{H}_B(H/k)$, or simply $\mathcal{H}_B(H)$, be the set of homotopy
equivalence classes of right faithfully flat $H$-Galois extensions of~$B$.
\end{Def}

\begin{Rems}\label{homotopyremarks}
(1)  Homotopy of $H$-Galois extensions can be formulated without changing
the ground ring from $k$ to~$k[t]$. Let $A_0$ and $A_1$ be right faithfully
flat $H$-Galois extensions of~$B$. Then the following are equivalent:
\begin{enumerate}
\item [(a)] $A_0 \sim A_1$.
\item [(b)] There is a right faithfully flat $H$-Galois extension $A$ of
$B[t]$ with ground ring~$k$ such that $t$ lies in the center of~$A$ and
$$A/(t-i) \cong A_i$$
as $H$-Galois extensions of $B$ for $i \in \{0,1\}$.
\end{enumerate}

The implication (a) $\Rightarrow$ (b) is clear. To prove the converse,
assume~(b).
Then $A$ is a $k[t]$-algebra, and the $H$-comodule structure map
$$A \to A \o H \cong A \o_{k[t]} H[t]$$
is $k[t]$-linear since $t$ is an $H$-coinvariant element.
Thus $B[t] \subset A$ is an $H[t]$-Galois extension, and (a)~follows.

(2) Our definition of homotopy equivalence extends the definition
introduced in~\cite{K} for central Hopf Galois extensions.
If $A_0$ and $A_1$ are right faithfully flat and central $H$-Galois
extensions of a commutative $k$-algebra $R$, we can view $A_0$ and $A_1$ as
$R \o H$-Galois extensions of the ground ring $R$ as explained above.
Then the following are equivalent:
\begin{enumerate}
\item [(a)] There is a faithfully flat and central $H$-Galois extension $A$
of~$R[t]$ with $[i]_*A \cong A_i$ for $i\in \{0,1\}$.
\item [(b)] There is a faithfully flat and central $H$-Galois extension $A$
of~$R[t]$ with $R_{[i]} \o_{R[t]} A \cong A_i$ for $i \in \{0,1\}$.
\item [(c)] $A_0 \sim A_1$ as $R \o H$-Galois extensions of the ground
ring~$R$.
\end{enumerate}

In (b), $R_{[i]}$ denotes $R$ as an $R[t]$-algebra via the $R$-algebra map
$[i] : R[t] \to R$ mapping $t$ onto $i$ for $i \in \{0,1\}$, and
$R_{[i]} \o_{R[t]} A$ is an algebra  by componentwise multiplication.

The equivalence of (a) and (b) follows from the bijectivity of the algebra map
$R_{[i]} \o_{R[t]} A \to k_{[i]} \o_{k[t]}A$ defined by $r \o a \mapsto 1
\o ra$,
and (b) is equivalent to (c) by Remark~\ref{homotopyremarks}~(1).

(3) Homotopy equivalence is different from isomorphism.
See~\cite[Section 4.5]{K} for examples of non-isomorphic homotopy
equivalent Hopf Galois extensions.

(4) If the Hopf algebra $H$ is finitely generated and projective as a module
over the ground ring~$k$, then by a theorem of Kreimer and Takeuchi~\cite{KT}
any $H$-Galois extension $A$ of $B$ is finitely generated projective
as a $B$-module. Thus $A$ defines an element $[A]$ in the $K$-theory
group $K_0(B)$ of~$B$.  
It follows from the definition that, if $A_0$ and $A_1$ are faithfully
flat $H$-Galois extensions of $B$ such that $A_0 \sim A_1$, then there is
$\xi \in K_0(B[t])$ such that for $i \in \{0,1\}$,
$$K_0([i])(\xi) = [A_i] \in K_0(B).$$
If $B$ is a regular ring, then by \cite[Chapter~XII, Theorem~3.1]{B},
$$K_0([0]) = K_0([1]) : K_0(B[t]) \to K_0(B). $$
(The latter also holds for certain non-regular rings~$B$, see~\cite{P}.)
In this case, $[A_0] = [A_1]$ in $K_0(B)$, which means that
the element we have constructed in~$K_0(B)$ is
invariant under homotopy equivalence.

If $k$ is a field and the antipode of $H$ is bijective, 
then any faithfully flat $H$-Galois extension $A$ of $B$ is projective over $B$ 
by a recent result in \cite[2.4.9]{Sch}.
Hence the above construction can sometimes be extended to cases when
$H$ is no longer finite-dimensional over~$k$ (now assumed to be a field).
For instance, let $H$ be a cosemisimple Hopf algebra with a decomposition
$H \cong \oplus_{i\in I}\, C_i$, where each $C_i$ is a finite-dimensional
coalgebra over~$k$. Then $A \cong \oplus_{i\in I}\, A \Box_H C_i$. 
Each summand $A \Box_H C_i$ of~$A$
is a finitely generated $B$-module by faithfully flat descent, 
and it is projective over $B$ since $A$ is $B$-projective. We then obtain a family
of elements of~$K_0(B)$ for which we can argue as above.
\end{Rems}

Next we show that $()_*$ and $()^*$ induce maps on the homotopy classes.

\begin{Prop}\label{homotopy}
Let $H$ be a Hopf algebra, and $A_0$, $A_1$ right faithfully flat
$H$-Galois extensions of~$B$ with $A_0 \sim A_1$.
\begin{enumerate}
\item If $\alpha : R \to S$ is a homomorphism of commutative $k$-algebras,
then $\alpha_*A_0 \sim \alpha_*A_1$.
\item If $\varphi : K \to H$ is a homomorphism of Hopf algebras and $K$
is flat over $k$, then $\varphi^*A_0 \sim \varphi^*A_1$.
\item If $\sigma : H \o H \to k$ is an invertible 2-cocycle,
then $(A_0)^{\sigma} \sim (A_1)^{\sigma}$.
\end{enumerate}
\end{Prop}

\pf
Let $A$ be a right faithfully flat $H[t]$-Galois extension of $B[t]$
with $[i]_*A \cong A_i$, $i \in \{0,1\}$.

(1) Extend $\alpha$ to the ring homomorphism $\alpha[t] : R[t] \to S[t]$
by mapping $t$ onto~$t$. Then $\alpha [i] = [i]\alpha[t]$.
By Proposition~\ref{baseextension}, $\alpha[t]_*A$ is a right faithfully flat
$H[t]$-Galois extension of~$B[t]$, and by Proposition~\ref{commute} (1),
$$[i]_*\alpha[t]_*A \cong ([i] \alpha[t])_*A \cong
(\alpha[i])_*A \cong \alpha_* [i]_*A \cong \alpha_*A_i  \text{ for } i\in
\{0,1\}.$$

(2) Let $i\in \{0,1\}$. We apply part (3) of Proposition~\ref{commute} to
the ring homomorphism $\tau = [i] : k[t] \to k$ and
to the $k[t]$-Hopf algebra map $\varphi[t] : K[t] \to H[t]$.
Then the Hopf algebra map $k_{[i]} \o_{k[t]} \varphi$ defined by
ground ring extension of $\varphi$ with respect to $[i]$
can be identified with~$\varphi$. Hence by Proposition~\ref{commute}~(3),
$$[i]_* \varphi[t]^*A \cong \varphi^* [i]_*A \cong \varphi^*A_i.$$

(3) By extension $\sigma$ defines an invertible 2-cocycle $\sigma[t]$ of
$H[t]$
over the ground ring~$k[t]$. By Proposition~\ref{twisting}, $A^{\sigma[t]}$
is a right faithfully flat $H[t]$-Galois extension of~$B[t]$, and
$$[i]_*(A^{\sigma[t]}) \cong ([i]_*A)^{\sigma}
\cong (A_i)^{\sigma} \text{ for } i\in \{0,1\}.$$
\epf

By Proposition~\ref{homotopy}, $\mathcal{H}_{R \o B}(R \o H)$ is a covariant
functor in commutative $k$-algebras~$R$, and $\mathcal{H}_B(H)$ is a
contravariant functor in $k$-flat Hopf algebras~$H$.

We now introduce homotopies between homomorphisms of Hopf algebras.

\begin{Def}
Let $K$, $H$ be Hopf algebras, and $\varphi : K \to H$ and $\psi : K \to H$
Hopf algebra homomorphisms.
We say that $\varphi$ and $\psi$ are {\em homotopic} ($\varphi \sim \psi$)
if there exists a $k[t]$-Hopf algebra homomorphism $\varPhi : K[t] \to H[t]$
with
$$[0](\varPhi(x)) = \varphi(x) \text{ and } [1](\varPhi(x)) = \psi(x)$$
for all $x \in K$. The map $\varPhi$ is called a {\em homotopy}
between $\varphi$ and~$\psi$.

The Hopf algebra map $\varphi : K \to H$ is a {\em homotopy equivalence}
if there exists a Hopf algebra homomorphism
$\varphi' : H \to K$ with $\varphi \varphi' \approx \id_H$ and
$\varphi' \varphi \approx \id_K$, where $\approx$ is the equivalence relation
generated by~$\sim$.
\end{Def}

\begin{Rems}\label{Hopfhomotopy}
(1) A homotopy $\varPhi$ between $\varphi$ and $\psi$ is given by a family
$\varPhi_n : K \to H$ ($n\geq 0$) of $k$-linear maps such that
for all $x \in K$, $\varPhi_n(x) \neq 0$ only for finitely many~$n$,
and for all $x,y \in K$ and all~$n \geq 0$,
\begin{enumerate}
\item [(a)] $\varPhi_n(xy) = \sum_{i+j = n} \varPhi_i(x) \varPhi_j(y)$,
\item [(b)]$\Delta(\varPhi_n(x))
= \sum_{i+j = n} \varPhi_i(x\sw1) \o \varPhi_j(x\sw2)$,
\item [(c)]$\varPhi_n(1) = \delta_{n0}$,
\item [(d)]$\varepsilon(\varPhi_n(x)) = \delta_{n0} \varepsilon(x)$,
\item [(e)]$\varPhi_0 = \varphi$,
\item [(f)]$\sum_{n\geq0}\varPhi_n(x) = \psi(x)$.
\end{enumerate}
The homotopy $\varPhi$ corresponding to the family $(\varPhi_n)$ is defined by
$$\varPhi(x) = \sum_{n\geq 0} \varPhi_n(x)t^n$$
for all $x \in K$.
Note that any family $(\varPhi_n)$ of $k$-linear maps with $\varPhi_n(x)
\neq 0$
only for finitely many~$n$ and (a)--(d) defines Hopf algebra homomorphisms
$\varPhi_0 : K \to H$ and $\sum_{n\geq0} \varPhi_n  : K \to H$.

(2) Let $H = \oplus_{n\geq 0} H(n)$ be an $\mathbb{N}$-graded Hopf algebra.
Then the inclusion $H(0) \subset H$ is a homotopy equivalence of Hopf algebras.
For the proof let $\iota : H(0) \to H$ be the inclusion and $\pi : H \to H(0)$
the projection. Both maps are Hopf algebra homomorphisms, and
$\pi \iota = \id_{H(0)}$.
We use the same homotopy as before for graded commutative algebras,
and define a $k[t]$-Hopf algebra homomorphism
$$\varPhi : H[t] \to H[t] \text{ with } \varphi(x)
= xt^n \text{ for all } x\in H(n), n\geq0.$$
Note that $\varPhi$ is a coalgebra map since for all $n \geq 0$,
$$\Delta(H(n)) \subset \oplus_{i+j=n} H(i) \o H(j) \text{ and }
\varepsilon(H(n)) = 0 \text{ if } n>0.$$
Then $\varPhi$ is a homotopy between $\iota \pi$ and~$\id_H$.
\end{Rems}

The first part of the next theorem generalizes~\cite[Proposition~2.3]{K}.

\begin{Thm}\label{mainhomotopy}
Let $B$ be an algebra and $H$  a Hopf algebra over~$k$.
\begin{enumerate}
\item If $\alpha : R \to S$ and $\beta : R \to S$ are homotopic
homomorphisms of
commutative $k$-algebras, then
$$\alpha_* = \beta_* : \mathcal{H}_{R \o B}(R \o H)
\to \mathcal{H}_{S \o B}(S \o H).$$
\item If $\varphi : K \to H$ and $\psi : K \to H$ are homotopic Hopf
algebra homomorphisms, where $K$ is a $k$-flat Hopf algebra, then
$$\varphi^* = \psi^* : \mathcal{H}_B(H) \to  \mathcal{H}_B(K).$$
\end{enumerate}
\end{Thm}

\pf
By Proposition~\ref{homotopy}, all maps are well defined on homotopy classes.

(1) Let $A$ be a right faithfully flat $R \o H$-Galois extension of~$R \o B$.
We will show that $\alpha_*A \sim \beta_*A$.

There is a homotopy $\theta : R \to S[t]$ with $[0]\theta = \alpha$ and
$[1]\theta =  \beta$.
Then $\theta_*A$ is a right faithfully flat $(S \o H)[t]$-Galois extension
of~$(S \o B)[t]$, and
$$[0]_*\theta_*A \cong ([0]\theta)_*A \cong \alpha_*A \text{ and }
[1]_*\theta_*A \cong ([1]\theta)_*A \cong \beta_*A.$$

(2) Let $A$ be a right faithfully flat $H$-Galois extension of $B$.
We will show that $\varphi^*A \sim \psi^*A$.

There is a homotopy $\varPhi : K[t] \to H[t]$ between $\varphi$ and $\psi$.
Define the $H[t]$-Galois extension $A[t]$ of $B[t]$ by ground ring extension
via~$k \subset k[t]$. Then $\varPhi^*A[t]$ is a right faithfully flat
$K[t]$-Galois  extension of~$B[t]$ by Proposition~\ref{cotensor},
and for $i \in \{0,1\}$ we have by Proposition~\ref{commute}
$$[i]_*\varPhi^*A[t] \cong (k_{[i]} \o_{k[t]} \varPhi)^*[i]_*A[t].$$
Since $A[t]$ is defined via ground ring extension with respect
to the inclusion $k \subset k[t]$,
we have $[i]_*A[t] \cong A$ for $i \in \{0,1\}$.
Let us identify the map $k_{[i]} \o_{k[t]} \varPhi$ defined
by ground ring extension via~$[i] : k[t] \to k$.
Note that for $i = 0,1$,
$$K \cong k_{[i]} \o_{k[t]} K[t], \;  x \mapsto 1 \o x$$
and
$$k_{[i]} \o_{k[t]} H[t] \cong H, \; 1 \o h t^n \mapsto hi^n$$
are isomorphisms. The image of an element $x \in K$ under the composition
$$K \cong k_{[i]} \o_{k[t]} K[t] \xrightarrow{k_{[i]} \o\varPhi}
k_{[i]} \o_{k[t]} H[t] \cong H$$
is $\sum_{n\geq0} \varPhi_n(x)i^n$,
where the homotopy is given by the family~$(\varPhi_n)_{n\geq 0}$.
Hence we can identify $k_{[0]} \o_{k[t]} \varPhi$  with $\varphi$,
and $k_{[1]} \o_{k[t]} \varPhi$ with~$\psi$.
This proves our claim that $\varphi^*A \sim \psi^*A$.
\epf

It is clear from Theorem \ref{mainhomotopy} that homotopy equivalences induce
bijective maps on homotopy classes. Hence Remarks~\ref{homotopicmaps}~(2)
and \ref{Hopfhomotopy}~(2) imply the following Corollary whose first part
generalizes \cite[Corollary~2.4]{K}.

\begin{Cor}\label{gradedhomotopy}
Let $B$ be an algebra and $H$ a Hopf algebra over $k$.
\begin{enumerate}
\item If $R = \oplus_{n\geq0} R(n)$ is an $\mathbb{N}$-graded commutative
algebra with $R(0) = k$, then the inclusion map $\iota : k \to R$ induces
a bijective map
$$\iota_* : \mathcal{H}_B(H) \to \mathcal{H}_{R \o B}(R \o H).$$
\label{gradedhomotopy1}
\item If $H = \oplus_{n\geq 0} H(n)$ is an $\mathbb{N}$-graded Hopf algebra
with $K=H(0)$, then the inclusion map $\iota : K \to H$ induces a bijective
map
$$\iota^* : \mathcal{H}_B(H) \to \mathcal{H}_B(K).$$\label{gradedhomotopy2}
\end{enumerate}
\end{Cor}

The following examples show that $\mathcal{H}$ cannot be replaced by $\Gal$
in Corollary~\ref{gradedhomotopy}~(1) and~(2).

\begin{Expls} (1) Let $H$ be a finite-dimensional Hopf algebra over a
field~$k$,
and $R$ a commutative graded $k$-algebra with $R(0) = k$ such that
there exists a non-cleft faithfully flat $R\o H$-Galois extension of $R$
over the ground ring~$R$. Since any $H$-Galois extension of $k$ is cleft,
the induced map $\Gal_k(H) \to \Gal_R(R \o H)$ is not bijective.
An example of this type is described in Proposition~\ref{explicitexample}
below with $H = k[C_2]$ and~$R = k[x^2,x^3]$.

(2) We recall Masuoka's computation of $\Gal_k(H)$ for the Taft Hopf algebra
in~\cite{Ma1}. Let $k$ be a field, $N>2$ be a natural number,
and $q$ a root of unity of order $N$ in~$k$.
The Taft Hopf algebra of dimension $N^2$ is the algebra
$$H = H_{N^2} = k \langle g, x \mid g^N = 1, \; x^N = 0, \; gxg^{-1} =
qx\rangle$$
with Hopf algebra structure defined by
$$\Delta(g) = g \o g \text{ and } \Delta(x) = 1 \o x + x \o g.$$
It is an $\mathbb{N}$-graded Hopf algebra with $H(n) = k[G]x^n$ ($n \geq 0$),
where $G = G(H)$ is the cyclic group of order $N$ generated by~$g$.
For any $r \in U(k)$ and $s \in k$ define the $H$-comodule algebra $A_{r,s}$ by
$$A_{r,s} = k\langle a,b \mid a^N = r, \; b^N=s, \; aba^{-1} = qb\rangle$$
with $H$-comodule algebra structure $\delta$ given by
$$\delta(a) = a \o g \text{ and } \delta(b) = 1 \o x + b \o g.$$

By \cite[Proposition 2.17, Lemma 2.19]{Ma1}, the map $(r,s) \mapsto A_{r,s}$
defines a bijection
$$U(k)/U(k)^N \times k \cong \Gal_k(H).$$

On the other hand, for any $r \in U(k)$ define the $k[G]$-comodule
algebra $A_{r}$ by
$$A_{r} = k\langle a \mid a^N = r\rangle \text{ with } \delta(a) = a \o g.$$
Then $r \mapsto A_{r}$ defines a bijection
$$U(k)/U(k)^N \cong \Gal_k(k[G]).$$
If $r \in U(k)$ and $s \in k$, then
$A_{r} \to A_{r,s}  \Box_{H} k[G], a \mapsto a \o g$,
is an isomorphism of $k[G]$-comodule algebras.
Hence the map $\Gal_k(H) \to \Gal_k(k[G])$ induced by
the inclusion $k[G] = H(0) \subset H$ can be identified with the projection
$U(k)/U(k)^N \times k \to U(k)/U(k)^N$. In particular, it is not bijective.

Without using the explicit computation of $\Gal_k(H)$, we know from
Corollary~\ref{gradedhomotopy}~(2) that the inclusion $k[G] \subset H$ defines
a bijection
$$\mathcal{H}_k(H) \cong \mathcal{H}_k(k[G]).$$
Since $G$ is cyclic of order $N$,
$$\mathcal{H}_k(k[G]) \cong \Gal_k(k[G]) \cong H^2(G,U(k)) \cong U(k)/U(k)^N$$
by Proposition~\ref{grouphomotopy}~(1) below.
\end{Expls}

For later use we note the following combination of our results
on twisted and graded Hopf algebras.

\begin{Cor}\label{gradedtwist}
Let $B$ be an algebra, $H = \oplus_{n\geq 0} H(n)$ an $\mathbb{N}$-graded
Hopf algebra, and~$K=H(0)$. Let $\sigma : H \times H \to k$ be an
invertible 2-cocycle such that
$\sigma(x,y) = \varepsilon(x)\varepsilon(y)$ for all $x,y \in K$.
Then $K$ is a Hopf subalgebra of $H^{\sigma}$ and
the inclusion $\iota : K \to H^{\sigma}$ induces a bijective map
$\iota^* : \mathcal{H}_B(H^{\sigma}) \to \mathcal{H}_B(K)$.
\end{Cor}

\pf
Let $A$ be a faithfully flat $H$-Galois extension of~$B$.
By Proposition~\ref{twisting} and Corollary~\ref{gradedhomotopy}~(2),
it is enough to show that
$$A \Box_H K \cong A^{\sigma} \Box_{H^{\sigma}} K$$
as $K$-Galois extensions of~$B$.
Since $A \Box_H K = A^{\sigma} \Box_{H^{\sigma}} K$ as $K$-comodules,
it suffices to check that the algebra structures of $A \Box_H K$ and
$A^{\sigma} \Box_{H^{\sigma}} K$ coincide.
Let $\sum_i a_i \o x_i$ and $\sum_j b_j \o y_j$ be elements of~$A \Box_H K$.
We denote the multiplication in $A^{\sigma} \Box_{H^{\sigma}} K$
by~$\cdot_{\sigma}$.
Then it follows from the definition of the cotensor product that
\begin{align}
\Bigl( \sum_i a_i \o x_i \Bigr) \cdot_{\sigma} \Bigl( \sum_j b_j \o y_j \Bigr)
&= \sum_{i,j} a_i \cdot_{\sigma} b_j \o x_i y_j  \notag \\
&= \sum_{i,j} a_i\sw0 b_j\sw0 \sigma^{-1}(a_i\sw1, b_j\sw1)  \o x_i y_j
\notag\\
&= \sum_{i,j} a_i b_j \sigma^{-1}(x_i\sw1, y_j\sw1)  \o x_i\sw2 y_j\sw2
\notag\\
&= \Bigl( \sum_i a_i \o x_i \Bigr) \Bigl( \sum_j b_j \o y_j \Bigr) \notag
\end{align}
since $\sigma^{-1}(x_i\sw1, y_j\sw1)=
\varepsilon(x_i\sw1)\varepsilon(y_j\sw1)$.
\epf

\section{Homotopy of abelian group functors, and examples from group algebras}

In \cite[Section 4]{Sw1} Swan constructed a homotopy functor out of any
functor
defined on a category of commutative algebras.
Our definition of $\mathcal{H}_B(H)$ is a special case of this construction.

A functor $F$ from the category of commutative $k$-algebras to the category of
sets is called a {\em homotopy functor} if for any commutative $k$-algebra
$R$,
$$F(R[0]) = F(R[1]) : F(R[t]) \to F(R)$$
or, equivalently, if the inclusion $\iota : R \subset R[t]$ induces a bijection
$F(R) \cong F(R[t])$.

For any commutative $k$-algebra $R$ let $\overline{F}(R)$ be the coequalizer
of the maps $[1]_*=F(R[1]),[0]_*= F(R[0]) : F(R[t]) \to F(R)$.

In general, if $M$, $N$ are sets and $f,g : M \to N$ are maps, the coequalizer
of the pair~$f,g$ is described as follows.
For all $x,y \in N$, define $x \sim y$ if there is an element $z \in M$ with
$f(z) = x$ and $g(z) = y$. Then the coequalizer of the pair $f,g$ is the
quotient map $N \to N/\!\!\approx$, where $\approx$ is the
equivalence relation generated by~$\sim$.

Let us say that two elements $x, y \in F(R)$ are {\em homotopy equivalent}
if~$x \approx y$.

From the definition there is a natural transformation $\eta : F \to
\overline{F}$.
By~\cite[Lemma~4.2]{Sw1}, $\eta$ is universal for maps of $F$ into homotopy
functors.
Thus $\overline{F}$ is the largest quotient of~$F$ which is a homotopy functor.

In particular, if $F$ is the functor
$$R \mapsto F(R)= \Gal_{R \o B}(R \o H),$$
then $\overline{F}(R) = \mathcal{H}_{R \o B}(R \o H)$.

\begin{Rem}\label{abelian}
Let $F$ be a functor from commutative $k$-algebras to abelian groups.
Let $\overline{F}$ be the largest quotient of the underlying set functor of
$F$
which is a homotopy functor. Then for any commutative algebra $R$,
the relation $\sim$ coincides with $\approx$ on~$F(R)$,
and $\overline{F}(R)$ is the cokernel of~$[1]_*- [0]_*$.

Indeed, using standard notation (see~\cite[Chapter~XII]{B}),
we denote by $NF(R)$ the kernel of the map $[0]_*: F(R[t]) \to F(R)$.
It splits off $F(R[t])$ and we have the following functorial decomposition:
$F(R[t]) = F(R) \oplus NF(R)$.
It is immediate to see that for two elements $x, y \in F(R)$
we have $x \sim y$ if and only if $x - y$ belongs to the image of the map
$$[1]_* : NF(R) \to F(R[t]) \to F(R),$$
which is equal to the image of $[1]_* - [0]_* : F(R[t]) \to F(R)$.
Therefore,
$$\overline{F}(R) \cong F(R)/[1]_*NF(R).$$
\end{Rem}

Let $H$ be a cocommutative $k$-flat Hopf algebra. Then the set~$\Gal_k(H)$
of isomorphism classes of Galois objects of $H$ form an abelian group
(see for example~\cite[10.5.3]{C}).
If $A$ and $A'$ are $H$-Galois objects, then their product in $\Gal_k(H)$
is the
isomorphism class of $A \Box _H A'$, where $A'$ is viewed as a left
$H$-comodule algebra, which is possible since $H$ is cocommutative.
In a very special case this group structure already appeared in~\cite{Ha}.
Thus
$$R \mapsto \Gal_R(R \o H)$$
is a functor from commutative $k$-algebras to abelian groups,
and $\mathcal{H}_k(H)$ is the cokernel of the homorphism
$$[1]_* - [0]_* : \Gal_{k[t]}(k[t] \o H) \to \Gal_k(H).$$

We are interested in the case when $H = R[G]$ is the group algebra of a
group $G$
over a commutative ring~$R$. Then Galois objects are $G$-strongly graded
algebras,
i.e., we have
$A = \oplus_{g \in G} A_g$, $A_1 = R$, and for all $g,h \in G$, the
multiplication map
$$A_g \o _R A_h \to A_gA_h = A_{gh}$$
is bijective (see~\cite[8.1.7]{M}).
Thus
$$\Phi(A) : G \to \Pic(R), \; g \mapsto [A_g],$$
is a group homomorphism.
Here $\Pic(R)$ is the abelian group of isomorphism classes of invertible
$R$-modules.

Suppose that for all $g \in G$, $A_g \cong R$ as an $R$-module. Then any
$A_g$ contains
an invertible element $u_g$ of $A$, and $A$ is $R[G]$-cleft with the
$R[G]$-colinear and invertible map $R[G] \to A$ defined by $g \mapsto u_g$.

Let $\sigma : G \times G \to U(R)$ be a 2-cocycle of the group $G$ acting
trivially
on~$U(R)$. The twisted group algebra  $\Psi(\sigma) = R_{\sigma}[G]$ is a
cleft $R[G]$-Galois object~\cite[7.1.5]{M}. Recall that $R_{\sigma}[G] = R[G]$
as a coalgebra with twisted multiplication $g \cdot_{\sigma} h =
\sigma(g,h) gh$
on basis elements $g,h \in G$.  The maps $\Phi$ and $\Psi$ define an exact
sequence of abelian groups
\begin{equation}\label{exact}
0 \to H^2(G,U(R)) \xrightarrow{\Psi} \Gal_R(R[G]) \xrightarrow{\Phi}
\Hom(G,\Pic(R)).
\end{equation}
Both homomorphisms $\Phi$ and $\Psi$ are natural transformations of abelian
group
functors on the category of commutative rings~$R$.

Let $\overline{\Hom}(G,\Pic(R))$ be the cokernel of
$$\Hom(G, [1]_* - [0]_*) : \Hom(G,\Pic(R[t])) \to \Hom(G,\Pic(R)).$$

\goodbreak

\begin{Prop}\label{grouphomotopy}
Let $G$ be a group and $R$ a reduced commutative ring.
\begin{enumerate}
\item If $\Pic(R[t]) = 0$, then
$$H^2(G,U(R)) \cong \Gal_R(R[G]) = \mathcal{H}_R(R[G]).$$
\item If $G$ is a finite abelian group, then $\Phi$ and $\Psi$ define an
exact sequence
$$H^2(G,U(R)) \to \mathcal{H}_R(R[G]) \to \overline{\Hom}(G,\Pic(R)) \to 0.$$
\end{enumerate}
\end{Prop}

\pf (1) The first isomorphism follows from \eqref{exact} since $\Pic(R) = 0$
as a direct summand in~$\Pic(R[t])=0$. It remains to show that homotopic
$R[G]$-Galois objects $A_0$ and $A_1$ are isomorphic. By assumption there is
an $R[t][G]$-Galois object~$A$  with $[i]_*A \cong A_i, i =0,1$.
Since $\Pic(R[t]) =0$, there is a 2-cocycle
$\sigma : G \times G \to U(R[t])$ with $A \cong R[t]_{\sigma}[G]$.
The ring $R$ being reduced, we have $U(R[t]) = U(R)$, and $\sigma$ takes
values
in~$U(R)$. Hence,
$$[i]_*A \cong [i]_*R[t]_{\sigma}[G] \cong R_{\sigma}[G]$$
for~$i = 0,1$.

(2) Define group homomorphisms $f_1, f_2,f_3$ by
$$f_1 = [1]_* - [0]_* :H^2(G,U(R[t])) \to H^2(G,U(R)),$$
$$f_2 = [1]_* - [0]_* : \Gal_{R[t]}(R[t][G]) \to \Gal_R(R[G]), \text{ and }$$
$$f_3 = \Hom(G,[1]_* - [0]_*) : \Hom(G,\Pic(R[t])) \to \Hom(G,R).$$
It is known (see~\cite[10.7.1]{C}) that the map $\Phi$ in~\eqref{exact} is an
epimorphism for any commutative ring~$R$ if the group $G$ is finite and
abelian.
Then $(f_1,f_2,f_3)$ define a homomorphism of short exact
sequences~\eqref{exact}
for $R[t]$ and for~$R$, and we have an exact sequence
$$\Coker(f_1) \to \Coker(f_2) \to \Coker(f_3) \to 0.$$
Since $R$ is reduced, $U(R) = U(R[t])$.
Therefore
$$[0]_* = [1]_* : U(R[t]) \to U(R),$$
hence $f_1 = 0$, and $\Coker(f_1) = H^2(G,U(R))$.
By Remark~\ref{abelian}, $\Coker(f_2) = \mathcal{H}_R(R[G])$,
and $\Coker(f_3) = \overline{\Hom}(G,\Pic(R))$.
\epf

Suppose in the situation of Proposition~\ref{grouphomotopy}~(2) that
$$[1]_* - [0]_* : \Pic(R[t]) \to \Pic(R)$$
is a split epimorphism. Then $\Psi$ induces an epimorphism
$$H^2(G,U(R)) \to \Gal_R(R[G]]) \to \mathcal{H}_R(R[G]).$$
Hence, if there is a non-cleft $R[G]$-Galois object,
then there is a non-cleft $R[G]$-Galois object which is {\em homotopically
trivial},
that is, homotopically equivalent to the group algebra.
We will now explicitly construct such an example.
\medskip
\goodbreak

Let $k$ be a field, $R$ a commutative $k$-algebra without zero divisors,
and $F$
the field of fractions of~$R$.
If $P$ is a fractional ideal of~$R$ (i.e., an $R$-submodule of~$F$)
such that~$P^2 = R$, then we can consider the $C_2$-graded $R$-submodule
$R*P$ of~$F[C_2]$ defined by
$$(R*P)_0 = R \text{ and } (R*P)_1 = Pg,$$
where $g$ generates the group~$C_2$.
The condition $P^2 = R$ implies that $R*P$ is a strongly $C_2$-graded
subalgebra
of~$F[C_2]$. It follows that $R*P$ is a $k[C_2]$-Hopf Galois extension
with coinvariants~$R$. The extension~$R*P$ is cleft if and only if $P$ is free.

We now assume that $k$ is of characteristic 2 and $R$ is not seminormal, i.e.,
there is an element $a \in F \setminus R$ such that $a^2$ and $a^3$ belong
to~$R$.
We also assume that any unit $u$ of the subring of~$F$ generated by $R$ and $a$
such that $u(1 + a) \in R$ belongs to~$R$.
As an example of such a ring, take $R$ to be the subring $k[x^2,x^3]$
in the polynomial algebra~$k[x]$, and~$a = x$.

Following~\cite[Section~7]{Sw2} and \cite[Section~2B, Example~2.15]{L},
we consider the Schanuel fractional ideal $P_a = (1 + a, a^2)$ of $R$
generated by $1 + a$ and $a^2$ in~$F$.
It follows from Equation~(2.15A) of~\cite{L} and the restriction we put
on the characteristic of~$k$ that~$P_a^2 = R$.
By the above we obtain a $k[C_2]$-Hopf Galois extension~$R*P_a$.
We claim the following.

\begin{Prop}\label{explicitexample}
The $k[C_2]$-Hopf Galois extension $R*P_a$ is not cleft, but it is
homotopically trivial.
\end{Prop}

\pf
(a) By~\cite[Proposition~2.15C]{L}, the $R$-module~$P_a$ is not free.
Therefore $R*P_a$ is not cleft.

(b) In the polynomial ring $R[t]$ the element $at$ of $F[t]$ is not in~$R[t]$,
but $a^2t^2$ and $a^3t^3$ belong to~$R[t]$.
We can then consider the Schanuel fractional ideal
$P_{at}$ of $R[t]$ generated by $1 + at$ and $a^2t^2$
in the field of fractions of~$R[t]$.
This leads to the non-cleft $k[C_2]$-Hopf Galois extension~$R[t]*P_{at}$
whose coinvariants are~$R[t]$.

We claim that $R[t]*P_{at} /(t-1) \cong R*P_a$ and $R[t]*P_{at} /(t) \cong
R[C_2]$.
This will prove that the Galois extension $R*P_a$ is homotopy equivalent
to the trivial extension~$R[C_2]$.
To prove the claim, it suffices to check that we have the following
isomorphisms of $R$-modules:
$$P_{at} /(t-1)P_{at} \cong P_a \text{ and } P_{at} /tP_{at} \cong R.$$
According to~\cite[Proof of Theorem~7.1]{Sw2},
the $R[t]$-module $P_{at}$ is the image of the idempotent matrix
$$M_t =
\begin{pmatrix}
 1 - a^4t^4 & a^2 t^2 + a^3 t^3 \\
 (1 + a^2 t^2)(a^2 t^2 - a^3 t^3) & a^4t^4
\end{pmatrix}
.$$
Therefore, $P_{at} /(t-1)P_{at}$ is the image of the idempotent
$$M_1 =
\begin{pmatrix}
 1 - a^4 & a^2  + a^3 \\
 (1 + a^2 )(a^2  - a^3) & a^4
\end{pmatrix}
,$$
which is the $R$-module $P_a$, and $P_{at} /tP_{at}$ is the image of
the idempotent
$M_0 =
 \begin{pmatrix}
 1 & 0 \\
 0 & 0
\end{pmatrix}$,
which is a free $R$-module of rank~$1$.
\epf

A fractional ideal as above is an invertible (rank one) projective~$R$-module.
Hence it represents an element $[P]$ of order~2 in the Picard group
$\Pic(R)$ of~$R$.
As observed above, the Hopf-Galois extension~$R*P$ is cleft if and only
if $[P] = 0$ in~$\Pic(R)$.

We have seen in Remark \ref{abelian} that homotopy equivalence in~$\Pic(R)$
is different from equality provided $[1]_*N\Pic(R)$ is non-zero. A necessary
condition for the non-vanishing of~$[1]_*N\Pic(R)$ is the non-vanishing
of~$N\Pic(R)$,
which by \cite[Theorem~1]{Sw2} is equivalent to the reduced quotient
$R_{\rm red}$ of $R$ not being seminormal.

The construction of the Schanuel module $P_{at}$ leading to the example of
the Proposition above, together with \cite[(2.23B) and Theorem~2.23]{L},
shows there exist injective maps
$$\iota : R/J \to \Pic(R) \text{ and } \iota_t : R/J \to N\Pic(R)$$
respectively given by
$$\iota(\overline{r}) = [(1+ra,a^2)] \text{ and }
\iota_t(\overline{r}) = [(1+rat,a^2t^2)],$$
where $J = \{ b\in R \, | \, ab \in R \}$ is the so-called conductor.
Moreover, we have $\iota  = [1]_* \circ \iota_t$.
This implies that $\Pic(R)/[1]_*N\Pic(R)$ is a quotient of the
cokernel $\Pic(R)/(R/J)$ of~$\iota$,
thus providing an ``upper bound'' to the set of homotopy classes in~$\Pic(R)$.

\begin{Expl}
If $R = k[x^2,x^3]$, then $\Pic(R) \cong R/J \cong k$, and the composition
$\iota$
$$R/J \xrightarrow{\iota_t} N\Pic(R) \subset
\Pic(R[t]) \xrightarrow{[1]_* - [0]_*} \Pic(R)$$
is an isomorphism~\cite[(2.23C)]{L}.
Hence $[1]_* - [0]_* : \Pic(R[t]) \to \Pic(R)$ is a split epimorphism.
In particular, all elements of the Picard group are homotopy equivalent to~$0$.
\end{Expl}

 \section{Examples from quantum groups}

From now on we assume that the ground ring $k$ is a field.
In this Section we apply our previous results to the Drinfeld-Jimbo quantum
groups $U_q(\mathfrak{g})$ and related finite-dimensional Hopf algebras.
It turns out that the computation of the homotopy classes for these Hopf
algebras
can be reduced to the case of abelian group algebras.

The general idea is to present these Hopf algebras as 2-cocycle twists of
graded Hopf algebras. Such a presentation was given for finite-dimensional
quantum groups in \cite{AS},~\cite{Di}. We adapt this approach to deal
with $U_q(\mathfrak{g})$ and prove some slightly more general results
which may be of independent interest.

We begin with the definition of the generalized quantum double
via twisting following~\cite{DT}.

Let $U,A$ be Hopf algebras, and $\tau : U \o A \to k$ a {\em skew-pairing},
that is a linear map such that for all $u,v \in U$ and $a,b \in A$,
\begin{align}
\tau(uv,a) &= \tau(u,a\sw1) \tau(v,a\sw2),\label{skew1} \\
\tau(u,ab) &= \tau(u\sw2,a) \tau(u\sw1,b),\label{skew2}\\
\tau(u,1) &= \varepsilon(u), \; \tau(1,a) = \varepsilon(a).
\end{align}

We assume that the antipode of $A$ is invertible. Then $\tau$ is invertible
with
respect to convolution,
and its inverse $\tau^{-1}$ is given for all $u \in U$, $a \in A$ by
\begin{equation}\label{tauinverse}
\tau^{-1}(u,a) = \tau(S(u),a) = \tau(u,S^{-1}(a)).
\end{equation}
Skew-pairings can be equivalently described as Hopf algebra homomorphisms
$\varphi : U \to (A^0)^{\cop}$, where $A^0$ is the dual Hopf algebra
(see~\cite[Section 9.1]{M}), and for any Hopf algebra $H$ we denote
by $H^{\cop}$ the algebra $H$ with the opposite comultiplication.
The connection between $\tau$ and $\varphi$ is given by
$$\tau(u,a) = \varphi(u)(a)$$
for all $u \in U$ and  $a \in A$.
We define the associated 2-cocycle $\sigma_{\tau} =\sigma$ on the
tensor product Hopf algebra $U \o A$ by
\begin{equation}\label{skewcocycle}
\sigma(u \o a, v \o b) = \varepsilon(u) \tau(v,a) \varepsilon(b)
\end{equation}
for all $u,v \in U$ and $a,b \in A$.
By definition, $(U \o A)^{\sigma}$ is the \emph{generalized quantum double}
of $U$, $A$ and~$\tau$.

The multiplication $\cdot_{\sigma}$ in $(U \o A)^{\sigma}$ is given by
\begin{equation}\label{taumultiplication}
(u \o a) \cdot_{\sigma} (v \o b)
= u \tau(v\sw1, a\sw1) v\sw2 \o a\sw2 \tau^{-1}(v\sw3,a\sw3)b
\end{equation}
for all $u,v \in U$ and $a,b \in A$.

\medskip

We now generalize some results on twisting from Didt's thesis~\cite{Di}.

We consider another skew-pairing $\widetilde{\tau}$ with
associated 2-cocycle~$\widetilde{\sigma}$.
Then $\rho = \widetilde{\sigma} \sigma^{-1} : (U \o A)^{\sigma} \o (U \o
A)^{\sigma} \to k$
is a 2-cocycle on the twisted Hopf algebra $(U \o A)^{\sigma}$ associated
to the
skew-pairing $\widetilde{\tau} \tau^{-1}$, and
$$(U \o A)^{\widetilde{\sigma}} = ((U \o A)^{\sigma})^{\rho}.$$

\begin{Lem}\label{groupliketwisting}
Let $y \in G(U), g \in G(A)$.
Then the following are equivalent:
\begin{enumerate}
\item For all $a \in A$, $u \in U$,
$$\tau(y,a) = \widetilde{\tau}(y,a) \text{ and }
\tau(u,g) = \widetilde{\tau}(u,g).$$
\item For all $a\in A$, $u \in U$,
$$\rho(u \o a, y \o g) = \varepsilon(u) \varepsilon(a) = \rho(y \o g, u \o
a).$$
\end{enumerate}
\end{Lem}
\pf
This is easily checked since for all $a\in A, u \in U$,
$$\rho(u \o a, y \o g) = \varepsilon(u)\widetilde{\tau}(y,a\sw1)
\tau^{-1}(y,a\sw2),$$
and
$$\rho(y \o g, u \o a) = \widetilde{\tau}(u\sw1,g)
\tau^{-1}(u\sw2,g)\varepsilon(a).$$
\epf

We will use the previous Lemma to define a 2-cocycle on a quotient Hopf
algebra
modulo a central subgroup algebra.

\begin{Lem}\label{quotientcocycle}
Let $H$ be a Hopf algebra, and $\rho : H \o H \to k$ an invertible 2-cocycle.
Let $G$ be a subgroup of $G(H)$, and assume that $G$ is central in $H$ and
that
$\rho(g,x) = \varepsilon(x) = \rho(x,g)$ for all $g \in G$ and $x \in H$.
Then $G$ is central in $H^{\rho}$, and
$\rho$ induces an invertible 2-cocycle $\overline{\rho}$ of the quotient
Hopf algebra~$H/(k[G])^+H$ such that
$$(H/(k[G])^+H)^{\overline{\rho}} = H^{\rho}/(k[G])^+H^{\rho}.$$
\end{Lem}

\pf In order to see that the map
$$\overline{\rho} : H/(k[G])^+H \o H/(k[G])^+H \to k, \;
\overline{x} \o \overline{y} \mapsto \rho(x,y),$$
is well defined, we have to show that
$$\rho(gx,y) = \rho(x,y) = \rho(x,yg)$$
for all $x,y \in H$, $g \in G$.
This follows from the assumption and the two cases of the 2-cocycle
condition~\eqref{cocycle} when $x=g$ and $z=g$. The rest of the Lemma
is then obvious since $(k[G])^+H = (k[G])^+H^{\rho}$ by~\eqref{twisted}.
\epf

In the sequel we will assume that the group $G(A)$ is abelian and that
there exist elements
$a_1, \ldots , a_t \in A \setminus 0$, $g_1, \ldots, g_t \in G(A)$,
and $\chi_1,\ldots,\chi_t \in \Hom(G(A), U(k))$
such that $A$ as an algebra is generated by $G(A)$ and
by~$a_1, \ldots , a_t$, and for all~$j$,
\begin{align}
\Delta(a_j) = g_j \o a_j + a_j \o 1 ,\label{A1} \\
ga_j g^{-1} = \chi_j(g) a_j \text{ for all } g \in G(A), \label{A2}\\
\chi_j(g_j) \neq 1.\label{A3}
\end{align}

Similarly, we assume that the group $G(U)$ is abelian and that there exist
elements $u_1, \ldots , u_s \in U \setminus 0$,  $y_1, \ldots, y_s \in G(U)$,
and $\eta_1,\ldots,\eta_s \in \Hom(G(U), U(k))$ such that $U$ as an algebra
is generated by $G(U)$ and by~$u_1, \ldots , u_s$, and for all~$i$,
\begin{align}
\Delta(u_i) = y_i \o u_i + u_i \o 1 , \label{U1}\\
yu_i y^{-1} = \eta_i(y) u_i \text{ for all } y \in G(U),  \label{U2}\\
\eta_i(y_i) \neq 1 .\label{U3}
\end{align}

Our assumptions for $U$ (and similarly for $A$) imply for all $i$
$$S(y_i) = y_i^{-1}, \; \varepsilon(y_i) = 1,
\;S(u_i) = -y_i^{-1} u_i, \;\varepsilon(u_i) = 0.$$

Note that, if $\tau : U \o A \to k$ is a skew-pairing, then
\begin{equation}\label{tauz}
\tau(y,a_j) = 0 \text{ for all } y \in G(U), 1 \leq j \leq t,
\end{equation}
since by \eqref{skew2} the map $\gamma : A \to k, a \mapsto \tau(y,a)$,
is an algebra map, and for any algebra map $\gamma$, $\gamma(a_j) = 0$
since $\chi_j(g_j) \neq 1$ by~\eqref{A3},  and by~\eqref{A2},
$\gamma(a_j) = \gamma(g_ja_jg_j^{-1}) = \chi_j(g_j) \gamma(a_j)$.

In the same way,
\begin{equation}\label{taug}
\tau(u_i,g) = 0 \text{ for all } g \in G(A), 1 \leq i \leq s.
\end{equation}

Generalizing an argument in the proof of~\cite[Theorem~5.17]{AS} on page~17,
we describe central group-like elements in~$(U \o A)^{\sigma}$.

\begin{Lem}\label{center}
Let $y \in G(U), g \in G(A)$.  Then the following are equivalent:
\begin{enumerate}
\item The group-like element $y \o g^{-1}$ is central in $(U \o A)^{\sigma}$.
\item For all  $1 \leq i \leq s$, $1 \leq j\leq t$,
$$\eta_i(y) = \tau(y_i,g) \text{ and } \chi_j(g^{-1}) = \tau(y,g_j).$$
\end{enumerate}
\end{Lem}

\pf
Since for all $u \in U$, $a \in A$,  $u \o a = (u \o 1) \cdot_{\sigma} (1
\o a)$,
the element $y \o g^{-1}$ is central if it commutes with all $u' \o 1$,
where $u'$ is group-like or $u' = u_i$ for some~$i$,
and with all $1 \o a'$, where $a'$ is group-like or $a' = a_j$ for some~$j$.

Since the groups $G(U)$ and $G(A)$ are abelian, $y \o g^{-1}$ commutes with
all $u' \o a'$, where $u'$ and $a'$ are group-like.

Using~\eqref{U1}, we compute
\begin{align}\notag
(y \o g^{-1}) \cdot_{\sigma} (u_i \o 1)
&= y \tau(u_{i}\sw1, g^{-1}) u_i\sw2 \o g^{-1}
\tau^{-1}(u_i\sw3,g^{-1})\\\notag
&= y \tau(y_i,g^{-1})y_i \o g^{-1} \tau^{-1}(u_i,g^{-1}) \\\notag
&\phantom{==}+y \tau(y_i,g^{-1}) u_i \o g^{-1} \tau^{-1}(1,g^{-1}) \\\notag
&\phantom{==}+y \tau(u_i,g^{-1}) \o g^{-1} \tau^{-1}(1,g^{-1})\\\notag
&=y\tau(y_i,g^{-1})u_i \o g^{-1},
\end{align}
since the first and last summand vanish by~\eqref{taug}. Since
$$(u_i \o 1)\cdot_{\sigma}(y \o g^{-1}) = u_i y\o g^{-1},$$
it follows from~\eqref{U2} that $y \o g^{-1}$ commutes with~$u_i \o 1$
if and only if $\eta_i(y) = \tau(y_i,g)$.
Similarly, $y \o g$ commutes with $1 \o a_j$ if and only
$\chi_j(g^{-1}) = \tau(y,g_j)$.
\epf

The next Proposition is an immediate consequence of the previous
Lemmas \ref{groupliketwisting}, \ref{quotientcocycle},~\ref{center},
together with~\eqref{tauz} and~\eqref{taug}.
It is formulated in~\cite{Di} in a special case for the Hopf algebras
studied in~\cite{AS}.

\begin{Prop}\label{generaltwistequivalent}
Let $U$, $A$ be Hopf algebras as above and
$$\tau, \widetilde{\tau} : U \o A \to k$$
invertible skew-pairings with associated 2-cocycles
$\sigma$,~$\widetilde{\sigma}$.
Assume that $s = t$, and
\begin{align}
\tau(y_i,g) = \widetilde{\tau}(y_i,g), \;
\tau(y, g_i) = \widetilde{\tau}(y,g_i),\label{quotientcondition}\\
\tau(y_i,g_j) = \chi_j(g_i^{-1})= \eta_i(y_j),
\end{align}
for all $y \in G(U)$, $g \in G(A)$, $1 \leq i,j \leq t$.

Then the subgroup $G$ generated by all $y_i \o g_i^{-1}$, $1 \leq i \leq t$,
is central in the twisted Hopf algebras $H = (U \o A)^{\sigma}$ and
$\widetilde{H} = (U \o A)^{\widetilde{\sigma}}$,
and the quotient Hopf algebra $\widetilde{H}/(k[G])^+\widetilde{H}$
is a 2-cocycle-twist of~$H/(k[G])^+H$.
\end{Prop}

In the situation of Proposition~\ref{generaltwistequivalent} we will denote
the images of
$$u_i \o 1, y_i \o 1, 1 \o a_j ,1 \o g_j, 1 \leq i,j \leq t,
\text{ in } H=(U \o A)^{\sigma}$$
respectively by $u_i$, $y_i$, $a_j$, $g_j$. Then a somewhat lenghty
calculation
using the multiplication rule~\eqref{taumultiplication} together
with \eqref{A1}, \eqref{U1}, \eqref{tauz},~\eqref{taug} and the equality
$$\tau^{-1}(u_j,a_i) = \tau(S(u_j),a_i) = - \tau(y_j^{-1},g_i)
\tau(u_j,a_i)$$
shows that for all $1\leq i,j \leq t$,
\begin{equation}\label{linkingrelation}
a_i u_j - \eta_j(y_i) u_j a_i = \tau(u_j,a_i) (1 - y_jg_i).
\end{equation}
By another abuse of language we use the same symbols $u_i, y_i, a_j, g_j$
for the images of these elements in the quotient Hopf algebra $H/(k[G])^+H$.
Then $y_j = g_j$ for all $j$ in $H/(k[G])^+H$, and
by multiplying with $g_j^{-1}$ we obtain from \eqref{linkingrelation}
\begin{equation}\label{linkingrelationprime}
a_i u_j' - u_j' a_i = \tau(u_j,a_i) (g_j^{-1} - g_i),
\end{equation}
with $u_j' = u_j g_j^{-1}$.

\medskip
\goodbreak

Let us now look at the concrete example of the Drinfeld-Jimbo
algebras $U_q(\mathfrak{g})$ (see~\cite[Chapter~4]{J}).
Let $k$ be a field of characteristic~$\neq 2,3$.
Let $(a_{ij})_{1 \leq i,j \leq t}$ be a Cartan
matrix of finite type of a semisimple complex Lie algebra $\mathfrak{g}$,
and for all $1 \leq i \leq t$ let $d_i \in \{1,2,3\}$ with $d_ia_{ij} = d_j
a_{ji}$
for all~$i,j$. Let $0 \neq q \in k$ and $q_i = q^{d_i} $ for all $i$.
Assume that $q^{2d_i} \neq 1$ for all~$i$.

Let $A = U_q^{\geq 0}(\mathfrak{g})$ be the algebra with generators
$E_i$, $K_i$, $K_i^{-1}$, $1 \leq i \leq t$ and relations
\begin{equation}
K_iK_j = K_j K_i, \;  K_i K_i^{-1} = 1 = K_i^{-1}K_i
\text{ for all } i,j,\label{grouprelations}
\end{equation}
\begin{align}
K_i E_j K_i^{-1} = q^{d_i a_{ij}} E_j
\text{ for all } i,j, \label{actionrelationsE}\\
\sum_{r =0}^{1 - a_{ij}}  (-1)^r \left[
\begin{matrix}
1 - a_{ij}\\
r
\end{matrix}
\right]_{q_i}
E_i^{1 - a_{ij} -r} E_j E_i^r &= 0 \text{ for all } i \neq j.\label{SerreE}
\end{align}
The algebra $A$ is a Hopf algebra with comultiplication
\begin{equation}\label{comultiplicationE}
\Delta(K_i) = K_i \o K_i, \; \Delta(E_i) = K_i \o E_i + E_i \o 1 \text{
for all }i.
\end{equation}

Let $U = U_q^{\leq 0}(\mathfrak{g})$ be the algebra with generators
$F_i$, $K_i$, $K_i^{-1}$, $1 \leq i \leq t$ and
relations~\eqref{grouprelations}
and
\begin{align}
K_i F_j K_i^{-1} = q^{-d_i a_{ij}} F_j \text{ for all } i,j,
\label{actionrelationsF}\\
\sum_{r =0}^{1 - a_{ij}}  (-1)^r \left[
\begin{matrix}
1 - a_{ij}\\
r
\end{matrix}
\right]_{q_i}
F_i^{1 - a_{ij} -r} F_j F_i^r &= 0 \text{ for all } i \neq j.\label{SerreF}
\end{align}
The algebra $U$ is a Hopf algebra with comultiplication
\begin{equation}\label{comultiplicationF}
\Delta(K_i) = K_i \o K_i, \; \Delta(F_i) = 1 \o F_i + F_i \o K_i^{-1}
\text{ for all }i.
\end{equation}
We define
$$a_i =E_i, \; g_i = K_i, \text{ and  } \chi_i(K_j)
= q^{d_i a_{ij}} \text{ for all }i,j.$$
Similarly, let
$$u_i = F_i K_i, \; y_i = K_i, \text{ and }\eta_i(K_j)
= q^{-d_i a_{ij}}\text{ for all }i,j.$$
Then all the conditions \eqref{A1}, \eqref{A2} and \eqref{A3} for~$A$
and \eqref{U1}, \eqref{U2} and \eqref{U3} for~$U$ are satisfied.

Let $\lambda_1, \ldots, \lambda_t$ be arbitrary elements in $k$, and
$\lambda = (\lambda_1, \ldots, \lambda_t)$. For all~$i$, define an algebra map
$$\gamma_i : A \to k \text{ with } \gamma_i(E_j) = 0 \text{ and }
\gamma_i(K_j) =\eta_i(K_j)  \text{ for all } j,$$
and an $(\varepsilon,\gamma_i)$-derivation
$$\delta_i^{\lambda} : A \to k \text{ with }
\delta_i^{\lambda}(E_j)= \delta_{ij} \lambda_i  \text{ and }
\delta_i^{\lambda}(K_j) = \delta_i^{\lambda}(K_j^{-1}) = 0 \text{ for all }
j.$$
Here, $\delta_{ij}$ is the Kronecker~$\delta$.
Finally we define a Hopf algebra homomorphism
$$\varphi^{\lambda} : U \to (A^0)^{\cop} \text{ by }
\varphi^{\lambda}(u_i)= \delta_i^{\lambda} \text{ and }
\varphi^{\lambda}(K_i) = \gamma_i\text{ for all } i.$$
The arguments in the proof of~\cite[Lemma~5.19]{AS} show the existence of
the algebra
maps~$\gamma_i$, the skew-derivations~$\delta_i^{\lambda}$
and the Hopf algebra map~$\varphi^{\lambda}$.

Let $\tau^{\lambda}$ and $\sigma^{\lambda}$ be the corresponding skew pairing
and 2-cocycle, and $H^{\lambda} = (U \o A)^{\sigma^{\lambda}}$.
The subgroup $G$ generated by all $K_i \o K_i^{-1}$ is central in
$H^{\lambda}$
by Lemma~\ref{center}. The quotient Hopf algebra
$$U_q^{\lambda}(\mathfrak{g}) = H^{\lambda}/(k[G])^+ H^{\lambda}$$
is generated as an algebra by elements $E_i$, $F_i$, $K_i$, $K_i^{-1}$,
$1\leq i \leq t$, with relations \eqref{grouprelations},
\eqref{actionrelationsE},
\eqref{SerreE}, \eqref{actionrelationsF}, \eqref{SerreF}, together with
\begin{equation}\label{linkingrelationslambda}
E_iF_j - F_jE_i = \delta_{ij} \lambda_i (K_i^{-1} - K_i) \text{ for all } i,j.
\end{equation}
This can be seen in the same way as in the end of the proof
of~\cite[Theorem 5.17]{AS}. Note that the relations
\eqref{linkingrelationslambda} are a special case of the
relations~\eqref{linkingrelationprime}.
The comultiplication is defined by
\eqref{comultiplicationE},~\eqref{comultiplicationF}.
The group $G(U_q^{\lambda}(\mathfrak{g}))$ is the free abelian group
with basis $K_1, \ldots, K_t$.

Since by definition, $U_q(\mathfrak{g}) = U_q^{\lambda}(\mathfrak{g})$,
where $\lambda_i = q^{-d_i} - q^{d_i}$ for all~$i$, we obtain the following.

\begin{Thm}\label{Uq}
Let $G = G(U_q(\mathfrak{g}))$. Then for any $k$-algebra $B$, the inclusion
$\iota : k[G] \to U_q(\mathfrak{g})$ induces a bijective map
$$\mathcal{H}_B(U_q(\mathfrak{g})) \cong \mathcal{H}_B(k[G]).$$
\end{Thm}

\pf
If we take $\lambda_i = 0$ for all~$i$, then $U_q^{\lambda}(\mathfrak{g})$
is a graded Hopf algebra. For another choice of~$\lambda$ we
get~$U_q(\mathfrak{g})$.
Hence by Proposition~\ref{generaltwistequivalent}, $U_q(\mathfrak{g})$ is
a 2-cocycle twist of a graded Hopf algebra with the group algebra $k[G]$
as degree 0 part, and the claim follows from Corollary~\ref{gradedtwist}.
\epf

\begin{Rem}\label{pointedingeneral}
Let $q$ be a primitive $N$-th root of unity, and assume that $N$ is odd
and not divisible by~3 if the Dynkin diagram of~$\mathfrak{g}$ contains
a component~$G_2$. The finite-dimensional Frobenius-Lusztig
kernel~$u_q(\mathfrak{g})$ is the quotient Hopf algebra of $U_q(\mathfrak{g})$
defined by adding the relations saying that the $N$-th powers of the positive
root vectors be 0 (see \cite[p.~120]{dCP},~\cite[p.~16]{AJS}).
By~\cite[Lemma~5.24]{AS} the above Hopf algebra homomorphism
$\varphi^{\lambda}$ also exists for $u_q(\mathfrak{g})$.
This implies that $u_q(\mathfrak{g})$ is a 2-cocycle twist of a graded
Hopf algebra. Hence we have a bijection
$$\mathcal{H}_B(u_q(\mathfrak{g})) \to \mathcal{H}_B(k[G]),$$
where $G = G(u_q(\mathfrak{g})) \cong (C_N)^t$.

There are numerous other Hopf algebras that are 2-cocycle twists of
graded Hopf algebras (see \cite{Ma2}, and~\cite{Di} for all the Hopf algebras
introduced in~\cite{AS}).  We do not know any example of pointed Hopf algebras
not of this type.
\end{Rem}

\end{document}